\documentclass{elsart}
\usepackage{ifpdf}
\usepackage{psfig}
\usepackage{graphicx,amssymb,lineno}
\ifpdf
\usepackage[%
  pdftitle={Instructions for use of the document class
    elsart},%
  pdfauthor={Simon Pepping},%
  pdfsubject={The preprint document class elsart},%
  pdfkeywords={instructions for use, elsart, document class},%
  pdfstartview=FitH,%
  bookmarks=true,%
  bookmarksopen=true,%
  breaklinks=true,%
  colorlinks=true,%
  linkcolor=blue,anchorcolor=blue,%
  citecolor=blue,filecolor=blue,%
  menucolor=blue,pagecolor=blue,%
  urlcolor=blue]{hyperref}
\else
\usepackage[%
  breaklinks=true,%
  colorlinks=true,%
  linkcolor=blue,anchorcolor=blue,%
  citecolor=blue,filecolor=blue,%
  menucolor=blue,pagecolor=blue,%
  urlcolor=blue]{hyperref}
\fi

\makeatletter
\def\elsartstyle{%
    \def\normalsize{\@setfontsize\normalsize\@xiipt{14.5}}
    \def\small{\@setfontsize\small\@xipt{13.6}}
    \let\footnotesize=\small
    \def\large{\@setfontsize\large\@xivpt{18}}
    \def\Large{\@setfontsize\Large\@xviipt{22}}
    \skip\@mpfootins = 18\p@ \@plus 2\p@
    \normalsize
}
\@ifundefined{square}{}{\let\Box\square}
\makeatother

\def\ds{\displaystyle}
\def\R{{\rm I\kern-.17em R}}
\begin{document}
\begin{frontmatter}
\title{A weighted finite difference method 
for the fractional diffusion equation based on the Riemann-Liouville derivative}

\author[label1]{Erc\'\i lia Sousa}
\author[label2]{Can Li \thanksref{label3}}
\address[label1]{CMUC, Department of Mathematics, University of Coimbra, 
3001-454 Coimbra, Portugal}
\address[label2]{School of Mathematics and statistics, Lanzhou University, Lanzhou 730000,P.R. China}
\thanks[label3]{Can Li was partly  supported by the Program for New Century Excellent
Talents in University under Grant No. NCET-09-0438, the National
Natural Science Foundation of China under Grant No. 10801067, and
the Fundamental Research Funds for the Central Universities under
Grant No. lzujbky-2010-63.}

\begin{abstract}
  A one dimensional fractional diffusion model with the
  Riemann-Liouville fractional derivative is studied.  First, a second
  order discretization for this derivative is presented and then an
  unconditionally stable weighted average finite difference method is
  derived.  The stability of this scheme is established by von Neumann
  analysis.  Some numerical results are shown, which demonstrate
  the efficiency and convergence of the method.  Additionally, some
  physical properties of this fractional diffusion system are
  simulated, which further confirm the effectiveness of our method.
\end{abstract}

\end{frontmatter}

\section{Introduction}

Recently, a large number of applied problems have been formulated on
fractional differential equations and consequently considerable attention has
been given to the solutions of those equations.  Fractional space derivatives
are used to model anomalous diffusion or dispersion, a phenomenon observed in
many problems.  There are some diffusion processes for which the Fick's second
law fails to describe the related transport behavior. This phenomenon is
called anomalous diffusion, which is characterized by the nonlinear growth of
the mean square displacement, of a diffusion particle over time.  The
anomalous diffusions differ according to the values of $\alpha$, where
$\alpha$ is the
order of the fractional derivative.  Some works providing an introduction to
fractional calculus related to diffusion problems are, for instance,
\cite{ben2001, gor1998, met2000, met22000, sch2001, zas2002}.  In this work we will be
interested in the anomalous diffusion, called supperdiffusion, for $1<\alpha
<2$ and experimental evidence of this type of diffusion
is already reported in several works \cite{ben2000a, hua2006, pac2000, zho2003}.

Fractional derivatives are non-local opposed to the local behaviour of
integer derivatives. Therefore, different challenges appear when we
try to derive numerical methods for this type of equations.  Numerical
approaches to different types of fractional diffusion models are
increasingly appearing in literature. We can found recent work on
numerical solutions for the fractional diffusion equation describing
superdiffusion \cite{erv2006,liu2004, lyn2003, mee2006, she2005,
  su2011, tad2006} and also for several transport equations including
this type of diffusion \cite{mee2004,sousa2009,zha2007}.
Some other works consider subdiffusion, which is represented by a time
fractional derivative of positive order and less than one
\cite{den2011, yus2005}.  However, the challenges for these equations
are different from the ones that arise when we consider a space
fractional derivative of order $1 \leq \alpha <2$.

Numerical methods, for models with superdiffusion, have been obtained
with mathematical techniques which do not necessarily consider a
second order discretization for the fractional derivative to achieve
second order accuracy.  In this work, we present a second order
approximation for the fra\-ctional Rie\-mann-Liouville derivative of
order $\alpha$, $1<\alpha < 2$. This approach uses some of the tools
described in \cite{die2004, li2009} and also applied in
\cite{sousa2010} to derive an approximation for the Caputo fractional
derivative defined in bounded domains.  Here, we consider the
Riemann-Liouville fractional derivative in an unbounded domain and its
discretization is represented by a series instead of a finite sum. We prove the order of consistency
of this discretization is second order.

A weighted average finite difference $\tau$-scheme is considered, for
$\tau \in [1/2,1]$, which includes the Crank-Nicolson method
($\tau=1/2$) and the back forward Euler method ($\tau=1$).  The consistency and
stability of the $\tau$-scheme are established and we
prove the $\tau$-scheme is unconditionally stable.  Also for $\tau=1/2$ we have second
order accuracy in time and space as expected.

Consider the one-dimensional fractional diffusion equation \cite{ben2000a, hua2006, tad2006}
\begin{equation}
\frac{\partial u}{\partial t}(x,t) = d(x) \frac{\partial^\alpha u}{\partial
  x^\alpha}(x,t)+p(x,t)
\label{fde}
\end{equation}
on the domain $x \in \R$, where $1 < \alpha \leq 2$ and $d(x) > 0$, subject to
the initial condition
\begin{equation}
u(x,0) = f(x), \ x \in \R
\end{equation}
and to the boundary condition
\begin{equation}
u(x,t)= 0 \quad \mbox{as} \quad |x| \rightarrow \infty.
\label{bc}
\end{equation}

The usual way of representing the fractional derivatives is by the
Riemann-Liouville formula.
The Riemann-Liouville fractional derivative of order
$\alpha$, for $x\in [a,b]$, $-\infty \leq a < b \leq \infty$, is defined by
\begin{equation}
\frac{\partial^\alpha u}{\partial x^\alpha}(x,t)  = 
\frac{1}{\Gamma(n-\alpha)}\frac{\partial^n}{\partial x^n}
\int_{a}^{x} {u(\xi,t)}{(x-\xi)^{n-\alpha-1}}d\xi , \quad (n-1 < \alpha < n)
\label{rl}
\end{equation}
where $\Gamma(\cdot)$ is the Gamma function and $n=[\alpha]+1$, with
$[\alpha]$ denoting the integer part of $\alpha$.

The function $u(x,t)$ under consideration, that is, which is solution of
(\ref{fde}), should be such that the corresponding integral (\ref{rl})
converges. If the function $u(x,t)$ vanishes at infinity, as assumed when we
impose the boundary condition (\ref{bc}), we have absolute convergence of such
integrals for a wide class of functions \cite{sam1993}. However, these
functions do not necessarily need to vanish at infinity and we can found under
which conditions these integrals converge in \cite{sam1993} (section 14.3).
There are very complete works about the fractional calculus
\cite{kil2006, mil1993, old2002, pod1999, sam1993}, where the theoretical 
properties of this type of derivative are studied in detail.

Another way to represent the
fractional derivative is by the Gr\"{u}nwald-Letnikov formula,
that is,
\begin{equation}
\frac{\partial^\alpha u}{\partial x^\alpha} (x,t) =  \lim_{\Delta x
  \rightarrow 0}
\frac{1}{\Delta x^\alpha}\sum_{k=0}^{\left[\frac{x-a}{\Delta x}\right]}
(-1)^k
\left(\begin{array}{c}
\alpha \\
k
\end{array}
\right)u(x-k\Delta x,t). \quad (\alpha > 0)\label{gl1}
\end{equation}
The Gr\"{u}nwald-Letnikov approximation is often used to numerically
approximate the Riemann-Liouville derivative and it was the first
algorithm to appear for approximating fractional derivatives
\cite{old2002, pod1999}.  However, this approximation has consistency of order
one and also very frequently numerical approximations based in this
formula originate unstable numerical methods and henceforth 
a shifted Gr\"{u}nwald-Letnikov formula is used \cite{tad2006,mee2004}.

The plan of the paper is as follows. In section 2 we derive a numerical
approximation for the Riemann-Liouville derivative. The full discretisation of the
fractional diffusion equation is given in section 3, where 
a weighted finite difference method in time is applied with the weight $\tau \in [1/2,1]$. 
In section 4 we prove the convergence of the numerical method by showing  consistency and
stability. In the fifth section we present numerical results which confirm
the theoretical results and in the last
section we give some conclusions.

\section{The numerical method}

In this section we present a numerical approximation for the Riemann-Liouville
derivative and also the numerical method that gives an approximate solution to
the fractional diffusion equation.

\subsection{Approximation of the Riemann-Liouville derivative}

Let us consider the Riemann-Liouville derivative \cite{old2002,pod1999}, that is,
\begin{equation}
 \frac{\partial^\alpha u}{\partial x^\alpha}(x,t)  = 
  \frac{1}{\Gamma(2-\alpha)}\frac{\partial ^2}{\partial x^2}
\int_{-\infty}^{x} {u(\xi,t)}{(x-\xi)^{1-\alpha}}d\xi, \quad 1 < \alpha < 2.
\label{rl2inf}
\end{equation}

We define the mesh points
$
x_j=j\Delta x, \ j \in \mathbb{Z} 
$
where $\Delta x$ denotes the  uniform space step.
For a fixed time $t$, let us denote
\begin{equation}
{\cal I}_\alpha(x) = \int_{-\infty}^{x} {u(\xi,t)}{(x-\xi)^{1-\alpha}}d\xi.
\label{ia}
\end{equation}

First, we do the following approximation at $x_j$ 
$$
\frac{\partial^2}{\partial x^2}{\cal I}_\alpha(x_j) \simeq \frac{1}{\Delta x^2}
\left[{\cal I}_\alpha(x_{j-1})-2{\cal I}_\alpha(x_{j})+{\cal I}_\alpha(x_{j+1})
\right].
$$

For each $x_j$ we need to calculate ${\cal I}_\alpha(x_j)$.

We compute these integrals by approximating $u(\xi,t)$, at a fixed instant $t$,
by a linear spline $s_j(\xi)$, whose nodes and knots are chosen at
$x_k$, $k=\dots, j-1, j$, that is, an approximation to 
${\cal I}_\alpha(x_j)$ becomes $ I_\alpha(x_j)$ defined by
\begin{equation}
I_\alpha(x_j)=\int_{-\infty}^{x_j}{s_j}(\xi){(x_j-\xi)^{1-\alpha}}d\xi.
\label{integraljinf}
\end{equation}
The spline $s_j(\xi)$ interpolates the points
$
\{(x_{k},t): \  k \leq j \} 
$
and is of the form \cite{pow1981}
\begin{equation}
s_j(\xi) = \sum_{k=-\infty}^{j} u(x_k,t) s_{j,k}(\xi),
\label{splineinf}
\end{equation}
with
$s_{j,k}(\xi)$, in each interval $[x_{k-1}, x_{k+1}]$, for $k \leq
j-1$, given by
\begin{equation}
s_{j,k}(\xi)=\left\{
\begin{array}{cc}
  \ds{\frac{\xi-x_{k-1}}{x_k-x_{k-1}}}, & x_{k-1} \leq \xi \leq x_k \\
& \\
\ds{\frac{x_{k+1}-\xi}{x_{k+1}-x_k}}, & x_{k} \leq \xi \leq x_{k+1} \\
& \\
0 & \mbox{otherwise,}
\end{array}
\right.
\label{sjk}
\end{equation}
and for $k=j$,
\begin{equation}
s_{j,j}(\xi)=\left\{
\begin{array}{cc}
\ds{\frac{\xi-x_{j-1}}{x_j-x_{j-1}}}, & x_{j-1} \leq \xi \leq x_j \\
& \\
0 & \mbox{otherwise.}
\end{array}
\right.
\label{sjj}
\end{equation}
From (\ref{integraljinf}) and (\ref{splineinf}),
\begin{equation}
I_\alpha(x_j)  =  \sum_{k=-\infty}^{j} u(x_k,t) \int_{x_{k-1}}^{x_{k+1}}
s_{j,k}(\xi)(x_j-\xi)^{1-\alpha} d\xi.
\label{ialphaxj0}
\end{equation}
We have that
\begin{eqnarray}
\int_{x_{k-1}}^{x_{k+1}}s_{j,k}(\xi)(x_j-\xi)^{1-\alpha} d\xi 
& = & \int_{x_{k-1}}^{x_{k}} \frac{\xi-x_{k-1}}{\Delta x}(x_j-\xi)^{1-\alpha}
+\int_{x_{k}}^{x_{k+1}}\frac{x_{k+1}-\xi}{\Delta x} (x_j-\xi)^{1-\alpha} {\nonumber} \\
& = & \frac{\Delta x^{2-\alpha}}{(2-\alpha)(3-\alpha)} a_{j,k},
\label{sjkxi}
\end{eqnarray}
where the $a_{j,k}$ are such that,
\begin{equation}
a_{j,k}=\left\{
\begin{array}{cc}
(j-k+1)^{3-\alpha}-2(j-k)^{3-\alpha}+(j-k-1)^{3-\alpha}, & \ k \leq j-1 \\
& \\
1, &  \ k=j.
\end{array}
\right.
\label{ajk}
\end{equation}
Therefore,
\begin{equation}
I_\alpha(x_j)  =   \frac{\Delta x^{2-\alpha}}{(2-\alpha)(3-\alpha)} 
\sum_{k=-\infty}^{j}u(x_k,t)a_{j,k},
\label{ialphaxj}
\end{equation}
and an approximation for $\ds{\frac{\partial^2}{\partial x^2}{\cal
    I}_\alpha(x_j) }$, is given by,
\begin{equation}
\frac{1}{\Delta x^2}
\left[I_\alpha(x_{j-1})-2I_\alpha(x_{j})+I_\alpha(x_{j+1})
\right]
\label{diffialpha}
\end{equation}
that is,
$$
\frac{\Delta x^{-\alpha}}{(2-\alpha)(3-\alpha)}
\left[\sum_{k=-\infty}^{j-1}u(x_k,t)a_{j-1,k}-2\sum_{k=-\infty}^{j}u(x_k,t)a_{j,k}+
\sum_{k=-\infty}^{j+1}u(x_k,t)a_{j+1,k}\right]. 
$$

Let us assume there are approximations ${\bf U}^n:=\{ U_{j}^{n} \}$ to the
values $u(x_j,t_n)$, where $t_n=n\Delta t, \ n\geq 0$ and $\Delta t$ is the
uniform time-step.  

We define the fractional operator as
\begin{equation}
\delta_\alpha U_j^{n} = \frac{1}{\Gamma(4-\alpha)}
\left\{\sum_{k=-\infty}^{j+1} q_{j,k}U_k^n \right\},
\label{fracoperinf}
\end{equation}
where 
\begin{eqnarray}
q_{j,k} & = & a_{j-1,k}-2a_{j,k}+a_{j+1,k}, \quad k \leq j-1{\nonumber} \\
q_{j,j} & = & -2a_{j,j}+a_{j+1,j} {\nonumber}\\
q_{j,j+1} & = & a_{j+1,j+1}. 
\label{qjk}
\end{eqnarray}
Therefore, an approximation of (\ref{rl2inf}),
for $t=t_n$, can be given by
$
\ds{\frac{\delta_\alpha U_j^n}{\Delta x^\alpha}}.
$

We can also write the fractional operator (\ref{fracoperinf}) as
\begin{equation}
\delta_\alpha U_j^{n} =
\frac{1}{\Gamma(4-\alpha)}
\sum_{m=-1}^{\infty}q_{j,j-m}U_{j-m}^n.
\label{intsumqjm}
\end{equation}

{\bf Remark}: Note that for $\alpha=1$ and $\alpha=2$ the coefficients
(\ref{qjk}) are such that $q_{j,k} = 0$, for $k < j-1$. For $\alpha=1$,
$q_{j,j-1} = -1$, $q_{j,j} = 0$, $q_{j,j+1} = 1$ and for $\alpha = 2$,
 $q_{j,j-1} = 1$, $q_{j,j} = -2$, $q_{j,j+1} = 1$.

{\bf Remark}: The series (\ref{intsumqjm}) converges absolutely for each $1 <
\alpha < 2$ and for every bounded function $u(x,t)$, for a fixed $t$.
This result is a straightforward consequence of some results given in section 3 about the
convergence of the series of the $q_{j,j-m}$.

In this section we have considered a linear spline to approximate the integral
representation of the Riemann-Liouville derivative with the purpose of
obtaining a second order approximation.  In the next section we describe the
full discretisation of the differential equation.

\subsection{Weighted average finite difference methods}

We discretize the spatial $\alpha$-order derivative following the steps of
the previous section. The discretization in time consists of the
weighted average discretization.

We consider the time discretization  $0 \leq
t_n \leq T$. Additionally, let $d_j=d(x_j)$,  $p_j^{n}=p(x_j,t_{n})$.
For the uniform space step $\Delta x$ and time step $\Delta t$, let
$$
\mu_j^\alpha=\frac{d_j\Delta t}{{\Delta x}^\alpha}.
$$

From equation (\ref{fde}) we can arrive
at the explicit Euler and implicit Euler numerical methods,
respectively
\begin{equation}\label{eulere}
\frac{U_j^{n+1}-U_j^n}{\Delta t}=\frac{d_j}{\Delta x^\alpha}\delta_{\alpha}
U^n_{j}+p^n_j ,~~~ 
\end{equation}
\begin{equation}\label{euleri}
\frac{U_j^{n+1}-U_j^n}{\Delta t}=\frac{d_j}{\Delta x^\alpha}\delta_{\alpha}
U^{n+1}_{j}+p^{n+1}_j ,~~~ 
\end{equation}

Let (\ref{eulere}) multiplies $(1-\tau)$ and (\ref{euleri})
multiplies $\tau$. We obtain the following weighted $\tau$-scheme
\begin{equation}\label{wscheme}
{U_j^{n+1}-U_j^n}=\mu_j^\alpha\bigg\{ (1-\tau)
\delta_{\alpha}U^n_{j}+\tau\delta_{\alpha} U^{n+1}_{j}\bigg\}
+\tau \Delta t p^{n+1}_j+(1-\tau)\Delta t p^n_j,
\end{equation}
where $\tau \in [1/2,1]$.

Note that for $\alpha = 2$, the operator (\ref{fracoperinf}) is the central
second order operator $\delta^2 U_j^{n}$, that is,
$$\delta_\alpha U_j^{n}=U_{j+1}^{n}-2U_{j}^{n}+U_{j-1}^{n}.
$$
We have the following numerical method
\begin{equation}
\left(1-\tau{\mu_j^\alpha}\delta_\alpha\right)U_j^{n+1} = 
\left(1+(1-\tau){\mu_j^\alpha}\delta_\alpha\right)U_j^{n}
+\Delta t p_j^{n+\tau},
\label{cn}
\end{equation}
where 
$$
 p_j^{n+\tau} = \tau p_j^n+(1-\tau)p_j^{n+1}.
$$

\section{Convergence of the numerical scheme}

In this section we prove the convergence of the numerical method
by showing it is consistent and von Neumann stable.
First, we start to study the consistency of the numerical method and lastly we
present the stability results. 

\subsection{Consistency}

In the beginning of this section, for the sake of clarity, we omit the variable
$t$ and  we denote the partial derivative of $u$ in $x$ of order $r$ by $u^{(r)}$.

\newtheorem{theorem1}{Lemma}
\begin{theorem1}
Let $u \in C^{(4)}(\R)$. For $\xi \in [x_{k-1}, x_k]$,
$$
u(\xi)-s_{j,k}(\xi) = -\frac{1}{r!}\sum_{r=2}^{3}u^{(r)}(\xi)l_{k,r}(\xi)
-\frac{1}{4!}u^{(4)}(\eta_k)l_{k,r}(\xi),
\ \eta_k \in [x_{k-1}, x_k],
$$
where
$$
|l_{k,r}(\xi)| \leq \Delta x^r.
$$
\end{theorem1}

{\bf Proof}:
For $\xi \in [x_{k-1}, x_k]$,
$$
u(\xi)-s_{j,k}(\xi) = u(\xi)-\frac{x_k-\xi}{\Delta x}u(x_{k-1})-\frac{\xi-x_{k-1}}{\Delta x}u(x_{k}).
$$
Using Taylor expansions, we obtain
\begin{eqnarray*}
u(\xi)-s_{j,k}(\xi) & = & -\frac{1}{r!}\sum_{r=2}^{3}u^{(r)}(\xi)l_{k,r}(\xi) -\frac{1}{4!}u^{(4)}(\eta_k)l_{k,r}(\xi),
\end{eqnarray*}
where $l_{k,r}(\xi)$ are functions which depend on $\Delta x$ and $x_k$, given
by
\begin{eqnarray}
l_{k,r}(\xi) & = & \frac{x_k-\xi}{\Delta x}(x_k-\xi-\Delta
x)^r-\frac{\xi-x_k+\Delta x}{\Delta x}(x_k-\xi)^r  \\
& = & (x_k-\xi)^r 
+ \sum_{r=0}^{p-1} 
\left(\begin{array}{c}
r \\
p
\end{array}
\right)
(x_k-\xi)^{p+1}(-1)^{r-p}\Delta x^{r-p-1}.
\end{eqnarray}
It is easy to conclude that $|l_{k,r}(\xi)| \leq \Delta x^r$, for $\xi \in
[x_{k-1}, x_k]$.
$\square$

\newtheorem{theorem20}[theorem1]{Theorem}
\begin{theorem20}
({Order of accuracy of the approximation for the fractional derivative}):
Let
$u \in C^{(4)}(\R)$ and such that $u^{(4)}(x) = 0$, for $x \leq a$, being $a$
a real constant. We have that
$$
\frac{\partial^\alpha u}{\partial x^\alpha}(x_j)
- \frac{\delta_\alpha u}{\Delta x^\alpha}(x_j) = \epsilon_1(x_j)+\epsilon_2(x_j),
$$
where
$$
|\epsilon_1(x_j)| \leq  C_1 \Delta x^2 \quad \quad |\epsilon_2(x_j)| \leq C_2 \Delta x^{2}, 
$$
and $C_1$ and $C_2$ are independent of $\Delta x$.
\end{theorem20}

{\bf Proof}: It is straightforward to prove that we have
\begin{eqnarray*}
\frac{\partial^\alpha u}{\partial x^\alpha} (x_j)& = &
\frac{1}{\Gamma(2-\alpha)}\frac{\partial^2}{\partial x^2} {\cal
  I}_\alpha(x_j)\\
& = & \frac{1}{\Gamma(2-\alpha)}\frac{1}{\Delta x^2}
\left[ {\cal I}_\alpha(x_{j-1}) -2{\cal I}_\alpha(x_j) + {\cal
    I}_\alpha(x_{j+1})\right] + \epsilon_1(x_j),
\end{eqnarray*}
where $\epsilon_1(x_j) = {\cal O}(\Delta x^2)$.

Let us define the error $E_S(x_j)$, such that,
$$
{\cal I}_\alpha(x_{j-1}) -2{\cal I}_\alpha(x_j) + {\cal
    I}_\alpha(x_{j+1}) = {I}_\alpha(x_{j-1}) -2{I}_\alpha(x_j) + {
    I}_\alpha(x_{j+1}) +E_S(x_j).
$$
We have
\begin{eqnarray*}
\frac{\partial^\alpha u}{\partial x^\alpha} (x_j)  & = & 
\frac{1}{\Gamma(2-\alpha)}\frac{1}{\Delta x^2}\left[
{\cal I}_\alpha(x_{j-1}) -2{\cal I}_\alpha(x_j) + {\cal
    I}_\alpha(x_{j+1})\right] \\
&&+\frac{1}{\Gamma(2-\alpha)}\frac{1}{\Delta x^2} E_S(x_j)+\epsilon_1(x_j),
\end{eqnarray*}
that is
$$
\frac{\partial^\alpha u}{\partial x^\alpha} (x_j) =
\frac{\delta_\alpha u}{\Delta x^\alpha}(x_j)+ \epsilon_1(x_j) +\epsilon_2(x_j),
$$
where
$$
\epsilon_2(x_j) = \frac{1}{\Gamma(2-\alpha)}\frac{1}{\Delta x^2} E_S(x_j).
$$
We are now going to compute the error $E_S(x_j)$.
We have
\begin{eqnarray*}
  E_S(x_j) & = & \sum_{k=-\infty}^{j-1} \
\int_{x_{k-1}}^{x_{k}}(u(\xi)-s_{j-1,k}(\xi))(x_{j-1}-\xi)^{1-\alpha} d\xi\\
  &&-2 \sum_{k=-\infty}^{j} \ \int_{x_{k-1}}^{x_{k}}(u(\xi)-s_{j,k}(\xi))(x_j-\xi)^{1-\alpha} d\xi\\
  &&+\sum_{k=-\infty}^{j+1} \ \int_{x_{k-1}}^{x_{k}}(u(\xi)-s_{j+1,k}(\xi))(x_{j+1}-\xi)^{1-\alpha} d\xi.
\end{eqnarray*}
Taking in consideration the previous lemma, let us denote
\begin{equation}
E_S(x_j) = -\sum_{r=2}^{4}\frac{1}{r!}E_r(x_j),
\label{es}
\end{equation}
where $E_r(x_j)$ are defined as follows.
For $r=2$ and $r=3$,
\begin{eqnarray}
  E_r(x_j) & = & \sum_{k=-\infty}^{j-1} \
\int_{x_{k-1}}^{x_{k}}l_{{k},r}(\xi)u^{(r)}(\xi)(x_{j-1}-\xi)^{1-\alpha} d\xi\nonumber\\
  &&-2 \sum_{k=-\infty}^{j} \ \int_{x_{k-1}}^{x_{k}}l_{{k},r}(\xi)u^{(r)}(\xi)(x_j-\xi)^{1-\alpha} d\xi\nonumber\\
  &&+\sum_{k=-\infty}^{j+1} \ \int_{x_{k-1}}^{x_{k}}l_{{k},r}(\xi)u^{(r)}(\xi)(x_{j+1}-\xi)^{1-\alpha} d\xi,
\label{e23}
\end{eqnarray}
and for $r=4$
\begin{eqnarray}
  E_r(x_j) & = & \sum_{k=-\infty}^{j-1} \
u^{(4)}(\eta_k)\int_{x_{k-1}}^{x_{k}}l_{{k},r}(\xi)(x_{j-1}-\xi)^{1-\alpha} d\xi\nonumber\\
  &&-2 \sum_{k=-\infty}^{j} u^{(4)}(\eta_k)\ \int_{x_{k-1}}^{x_{k}}l_{{k},r}(\xi)(x_j-\xi)^{1-\alpha} d\xi\nonumber\\
  &&+\sum_{k=-\infty}^{j+1} u^{(4)}(\eta_k)\ \int_{x_{k-1}}^{x_{k}}l_{{k},r}(\xi)(x_{j+1}-\xi)^{1-\alpha} d\xi.
\label{e4}
\end{eqnarray}
For $r=2,3$ by changing variables, we obtain
\begin{eqnarray*}
  E_r(x_j) & = & \sum_{k=-\infty}^{j} \
\int_{x_{k-1}}^{x_{k}}l_{k,r}(\xi)u^{(r)}(\xi-\Delta x)(x_{j}-\xi)^{1-\alpha} d\xi\\
  &&-2 \sum_{k=-\infty}^{j} \ \int_{x_{k-1}}^{x_{k}}l_{k,r}(\xi)u^{(r)}(\xi)(x_j-\xi)^{1-\alpha} d\xi\\
  &&+\sum_{k=-\infty}^{j} \
  \int_{x_{k-1}}^{x_{k}}l_{k,r}(\xi)u^{(r)}(\xi+\Delta x)(x_{j}-\xi)^{1-\alpha} d\xi,
\end{eqnarray*}
that is,
\begin{eqnarray*}
  E_r(x_j) & = &\sum_{k=-\infty}^{j} \ 
\int_{x_{k-1}}^{x_{k}}l_{k,r}(\xi)\left[
u^{(r)}(\xi+\Delta x)-2u^{(r)}(\xi)+u^{(r)}(\xi-\Delta x)\right](x_{j}-\xi)^{1-\alpha} d\xi.
\end{eqnarray*}

Let $x_a= N_a \Delta x$ such that $u^{(4)}(x)=0$, for $x \leq x_a$. 
For  $r=2$ we have
\begin{eqnarray*}
  E_2(x_j) & = &\sum_{k=-\infty}^{j} \ 
\int_{x_{k-1}}^{x_{k}}l_{k,2}(\xi)\left[
u^{(r)}(\xi+\Delta x)-2u^{(r)}(\xi)+u^{(r)}(\xi-\Delta
x)\right](x_{j}-\xi)^{1-\alpha} d\xi\\
& = & \frac{\Delta x^2}{2} \sum_{k=N_a+1}^{j}
u^{(4)}(\xi_k)c_{j,k,2}, \quad \xi_k \in [x_{k-1}, x_k]
\end{eqnarray*}
where
$$
c_{j,k,2}= \int_{x_{k-1}}^{x_{k}}l_{k,r}(\xi)(x_{j}-\xi)^{1-\alpha} d\xi
$$
Since, by Lemma 1,
$$
|c_{j,k,2}| \leq \Delta x^2 \int_{x_{k-1}}^{x_{k}}(x_{j}-\xi)^{1-\alpha} d\xi
$$
and
$$
\int_{x_{a}}^{x_{j}}(x_{j}-\xi)^{1-\alpha} d\xi = \frac{1}{2-\alpha}{(x_j-x_a)^{2-\alpha}}
$$
we have
\begin{equation}
|E_2(x_j)| \leq \frac{\Delta x^4}{2(2-\alpha)}||u^{(4)}||_\infty (x_j-x_a)^{2-\alpha}.
\label{de2}
\end{equation}

For $r=3$,
$$
E_3(x_j) = \sum_{k=N_a+1}^{j} \Delta
x(u^{(4)}(\xi_{k_1})-u^{(4)}(\xi_{k_2}))c_{j,k,3}, \ \xi_{k_1},\xi_{k_2} \in
[x_{k-1}, x_k]
$$
and
$$
|c_{j,k,3}|  \leq \Delta x^3 \int_{x_{k-1}}^{x_k} (x_j-\xi)^{1-\alpha} d\xi.
$$
We have
\begin{equation}
|E_3(x_j)| \leq  \frac{2\Delta x^4}{(2-\alpha)}||u^{(4)}||_\infty (x_j-x_a)^{2-\alpha}.
\label{de3}
\end{equation}
Finally for $r=4$, we bound each integral of (\ref{e4}) separately.
For the first integral we have
\begin{eqnarray*}
&&\sum_{k=N_a+1}^{j-1}u^{(4)}(\eta_k)\int_{x_{k-1}}^{x_k} l_{k,4}(\xi)(x_{j-1}-\xi)^{1-\alpha} d \xi\\
& \leq & 
\Delta x^4 ||u^{(4)}||_\infty \sum_{k=N_a+1}^{j-1}\int_{x_{k-1}}^{x_k}
(x_{j-1}-\xi)^{1-\alpha} d\xi \\
& = & \frac{\Delta x^4}{2-\alpha} ||u^{(4)}||_\infty(x_{j-1}-x_a)^{2-\alpha}.
\end{eqnarray*}
Therefore, since $(a+b)^p \leq |a|^p+|b|^p$ for $0<p\leq1$, we have
$$
\sum_{k=N_a+1}^{j-1}u^{(4)}(\eta_k)\int_{x_{k-1}}^{x_k} l_{k,4}(\xi)(x_{j-1}-\xi)^{1-\alpha} d \xi
\leq \frac{\Delta x^4}{2-\alpha}||u^{(4)}||_\infty ((x_{j}-x_a)^{2-\alpha}+\Delta x^{2-\alpha}).
$$
Similarly, for the second integral we have
$$
\sum_{k=N_a+1}^{j}u^{(4)}(\eta_k)\int_{x_{k-1}}^{x_k} l_{k,4}(\xi)(x_{j}-\xi)^{1-\alpha} d \xi
 \leq  \frac{\Delta x^4}{2-\alpha}
||u^{(4)}||_\infty(x_{j}-x_a)^{2-\alpha}
$$
and for the third integral
$$\sum_{k=N_a+1}^{j+1}u^{(4)}(\eta_k)\int_{x_{k-1}}^{x_k}
l_{k,4}(\xi)(x_{j+1}-\xi)^{1-\alpha} d \xi
\leq
\frac{\Delta x^4}{2-\alpha}||u^{(4)}||_\infty ((x_{j}-x_a)^{2-\alpha}+\Delta x^{2-\alpha}).
$$
Finally, we have
\begin{equation}
|E_4(x_j)| \leq \frac{3\Delta
  x^4}{2-\alpha}||u^{(4)}||_\infty (x_{j}-x_a)^{2-\alpha}+
\frac{2\Delta
  x^{6-\alpha}}{2-\alpha}||u^{(4)}||_\infty.
\label{de4}
\end{equation}
From (\ref{de2}), (\ref{de3}) and  (\ref{de4}) it is easy to conclude that the
error $E_S(x_j)$ defined by (\ref{es}) is of order $\mathcal{O}(\Delta x^4)$
and therefore the $\epsilon_2(x_j)$ is of order $\mathcal{O}(\Delta x^2)$.

$\square$

\newtheorem{theorem2}[theorem1]{Theorem}
\begin{theorem2}
The truncation error of the weighted numerical method (\ref{cn}) is of order $ {\cal O}(\Delta
x^2)+ {\cal O}(\Delta t^{m_\tau})$,
where
$m_\tau =1$, for $\tau \in (1/2,1]$ and $m_\tau=2$, for $\tau=1/2$.
\end{theorem2}
{\bf Proof}: Let $u=u(x,t)$ be a solution to the fractional partial
differential equation and satisfying the conditions of the previous theorem.
Note that the truncation error for the numerical
method (\ref{cn})  is given by
$$
T_j^n = \frac{u_j^{n+1}-u_j^n}{\Delta t} - 
\frac{d_j}{\Delta x^\alpha}\left(\tau\delta_\alpha u_j^{n+1}+(1-\tau)\delta_\alpha
  u_j^{n}\right)
-p_j^{n+\tau}.
$$
We have that
\begin{equation}\label{pro3.3}
\frac{u_j^{n+1}-u_j^n}{\Delta t}=\frac{\partial
u(x_{j},t_{n})}{\partial t}+\frac{\Delta t}{2}\frac{\partial^2
u(x_{j},t_{n})}{\partial t^2}+O(\Delta t^2),
\end{equation}
and using the previous theorem we have
\begin{eqnarray*}
T_j^n  & = & \frac{\partial u(x_{j},t_{n})}{\partial
t}+\frac{\Delta t}{2}\frac{\partial^2 u(x_{j},t_{n})}{\partial
t^2}+O(\Delta t^2)
-\tau\bigg(d_j\frac{\partial^{\alpha}u(x_{j},t_{n+1})}{\partial x^\alpha}
+O(\Delta x^2)\bigg)\\
&&-(1-\tau)\bigg(d_j\frac{\partial^{\alpha}u(x_{j},t_{n})}{\partial
x^\alpha}+O(\Delta x^2)\bigg)-p^{n+\tau}_j.
\end{eqnarray*}
Therefore
\begin{eqnarray*}
T_j^n &=&\frac{\partial u(x_{j},t_{n})}{\partial t}+\frac{\Delta
t}{2}\frac{\partial^2 u(x_{j},t_{n})}{\partial t^2}+-(1-\tau)\frac{\partial u(x_{j},t_{n})}{\partial
t}-\tau\frac{\partial u(x_{j},t_{n+1})}{\partial t}\\
&& +O(\Delta t^2) +O(\Delta x^2)
\end{eqnarray*}
Finally,
\begin{eqnarray*}
T_j^n &= & (\frac{1}{2}-\tau)\Delta t\frac{\partial^2
u(x_{j},t_n)}{\partial t^2}+O(\Delta t^2)+O(\Delta x^2).
\end{eqnarray*}

$\square$

\subsection{Fourier decomposition of the error}

In order to derive stability conditions for the finite difference schemes, we
apply the von Neumann analysis or Fourier analysis. Fourier analysis assumes
that we have a solution defined in the whole real line. It is also applied to
problems defined in  finite domains with periodic boundary conditions since 
the solution is seen as a
periodic function in $\R$.

If $u_j^n$ is the exact solution $u(x_j,t_n)$, let
\begin{equation}
E_j^n = U_j^n-u_j^n
\label{errorfourier}
\end{equation}
be the error at time level $n$ in mesh point $j$. To apply the von Neumann
analysis we also consider $d_j$ locally constant, and we denote $\mu_j^\alpha$
by $\mu^\alpha$.

Considering the scheme (\ref{cn}) and inserting equation
(\ref{errorfourier}) into that equation leads to
\begin{equation}
\left(1-\tau\mu^\alpha\delta_\alpha\right)E_j^{n+1} = 
\left(1+(1-\tau)\mu^\alpha\delta_\alpha\right)E_j^n .
\label{errorcn}
\end{equation}

The von Neumann analysis assumes that any finite mesh function, such as, the error
$E_j^n$ will be decomposed into a Fourier series as
$$
E_j^n = \sum_{p=-N}^{N} \kappa_p^n e^{i\xi_p(j\Delta x)}, \quad j=-N, \dots, N,
$$
where $\kappa_p^n$ is the amplitude of the $p$-th harmonic and
$\xi_p=p\pi/N\Delta x$. The product $\xi_p \Delta x$ is often called the phase
angle $\theta=\xi_p \Delta x$ and covers the domain $[-\pi,\pi]$ in steps of $\pi/N$.

Considering a single mode $\kappa^n e^{ij\theta}$, its time evolution is determined by
the same numerical scheme as the error $E_j^n$.  Hence
inserting a representation of this form into a numerical scheme we obtain
stability conditions. The stability conditions will be satisfied if the
amplitude factor $\kappa$ does not grow in time, that is, if we have
$|\kappa(\theta)| \leq 1, \ $ for all $\theta$.

As we have seen the fractional operator can be written as
$$
\delta_\alpha E_j^n = \frac{1}{\Gamma(4-\alpha)} \sum_{m=-1}^{\infty} q_{j,j-m}E_{j-m}^n,
$$
where the $q_{j,j-m}$ are defined by (\ref{qjk}).

\begin{figure}[h]
\centerline{
\psfig{figure=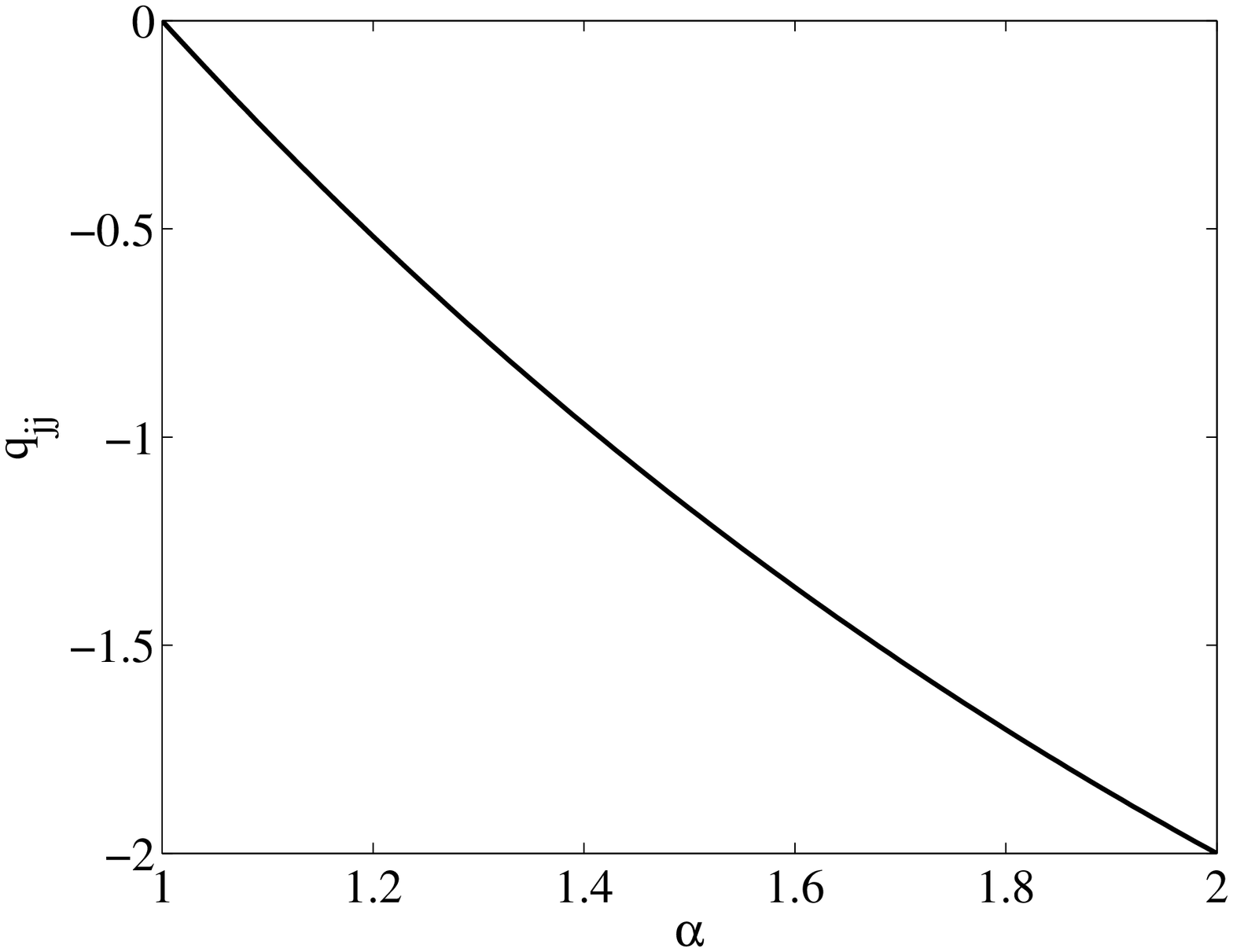,height=6.5cm,width=6.5cm}\quad \quad
\psfig{figure=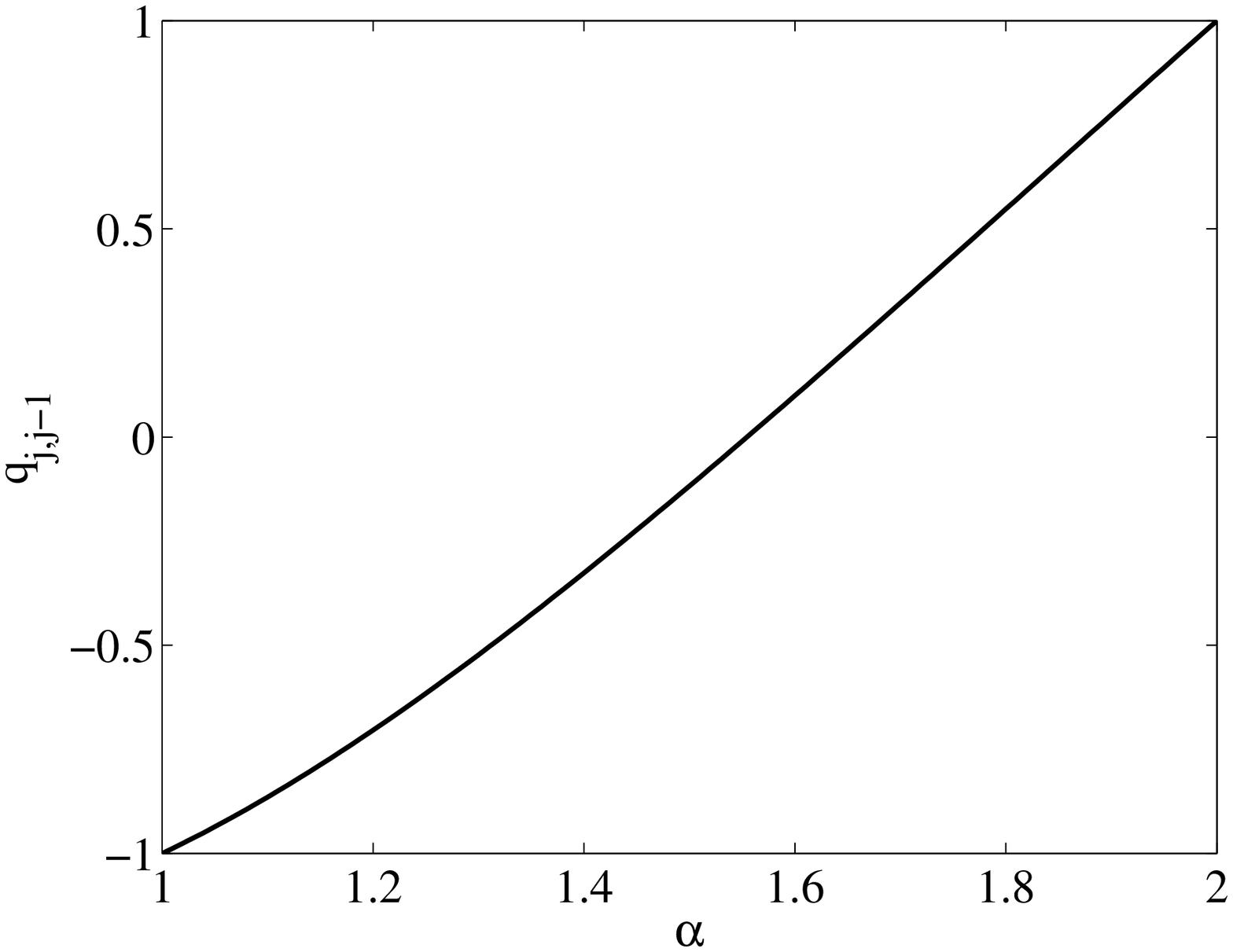,height=6.5cm,width=6.5cm}
}
\caption{\small Coefficients (\ref{qjk}): (a) $q_{jj}$
(b) $q_{j,j-1}$}
\label{d1}
\end{figure}

\begin{figure}[h]
\centerline{
\psfig{figure=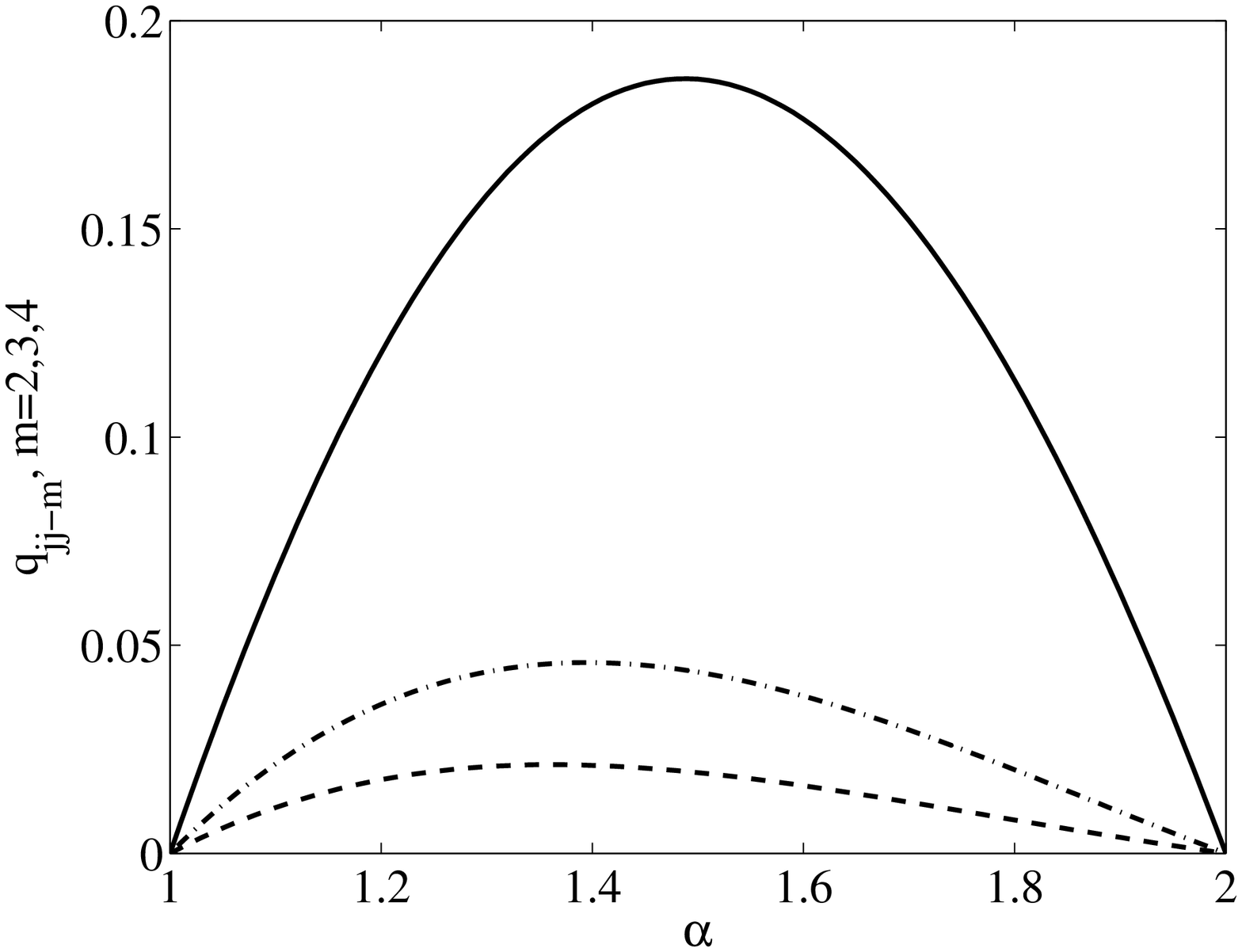,height=6.5cm,width=6.5cm}\quad \quad
\psfig{figure=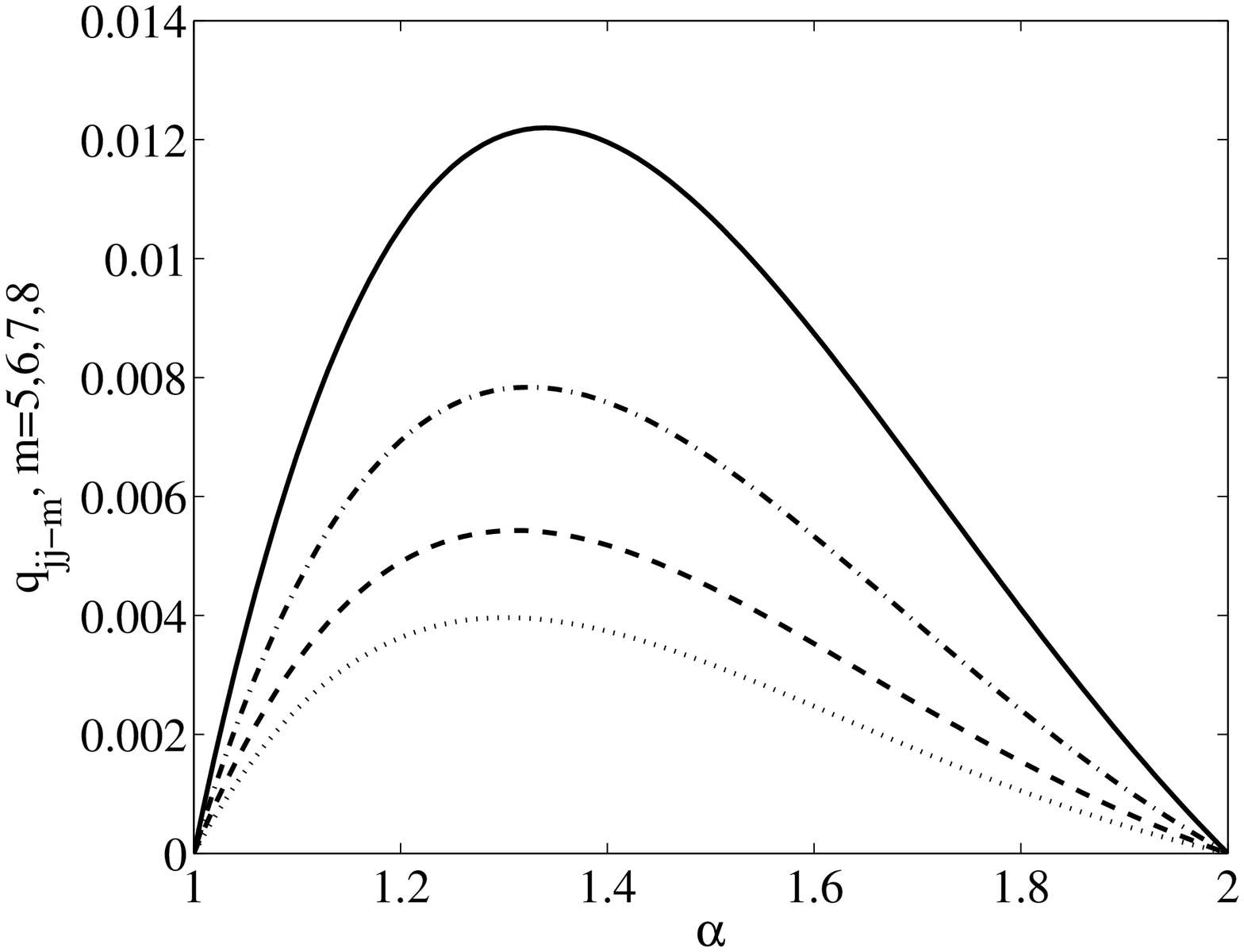,height=6.5cm,width=6.5cm}
}
\caption{\small Coefficients (\ref{qjk}): (a) $q_{j,j-m}$, $m=2,3,4$;
(b) $q_{j,j-m}$, $m=5,6,7,8$}
\label{d2}
\end{figure}

First we plot, in Figures \ref{d1} -- \ref{d2}, the coefficients $q_{j,j-m}$
and then we give the properties that
allow us to conclude this is a well-defined operator.

The following lemma characterizes the coefficients
$q_{j,j-m}$ and is useful to prove our next results.

\

\newtheorem{theorem4}[theorem1]{Lemma}
\begin{theorem4}
Consider the coefficients $q_{j,j-m}$ defined by (\ref{qjk}). Then

(a) $q_{j,j+1} = 1, \quad $ $ q_{j,j} \leq 0, \quad$ $ q_{j,j-m} \geq 0, \ m
\geq 2, \quad$ $ \ds{\lim_{m \rightarrow \infty} }q_{j,j-m} = 0 \quad $ and $q_{j,j-(m+1)} \leq
q_{j,j-m} \leq q_{j,j-2}$.

(b) $\ds{\sum_{m=2}^{\infty}} q_{j,j-m}= -3+3\times 2^{3-\alpha}-3^{3-\alpha}$.

(c) $\ds{\sum_{m=-1}^{\infty}} q_{j,j-m} = 0$.
\end{theorem4}

{\bf Proof} : (a) We have that 
$q_{j,j+1} = a_{j,j} =1$, $q_{j,j} = 2^{3-\alpha}-4 \leq 0$, for $ 1 <
\alpha \leq 2$ and $q_{j,j-1} = 3^{3-\alpha}-4 \times 2^{3-\alpha}+6$, which
can be positive or negative depending on the value of $\alpha$.
The $q_{j,j-m}$, $m\geq 2$, are of the form
$$
q_{j,j-m} =
(m+2)^{3-\alpha}-4(m+1)^{3-\alpha}+6m^{3-\alpha}-4(m-1)^{3-\alpha}+
(m-2)^{3-\alpha}.
$$
Hence,
\begin{eqnarray}
&&q_{j,j-m}  =  {\nonumber}\\
& & m^{3-\alpha}
\left[\left(1+\frac{2}{m}\right)^{3-\alpha}-4\left(1+\frac{1}{m}\right)^{3-\alpha}
+6-4\left(1-\frac{1}{m}\right)^{3-\alpha}+
\left(1-\frac{2}{m}\right)^{3-\alpha}
\right]{\nonumber}\\
&= &m^{3-\alpha}
\left[\sum_{k=0}^{\infty}\left(\begin{array}{c}
3-\alpha \\
k
\end{array}
\right)\left(\frac{2}{m}\right)^k
-4\sum_{k=0}^{\infty}\left(\begin{array}{c}
3-\alpha \\
k
\end{array}
\right)\left(\frac{1}{m}\right)^k+6\right.{\nonumber}\\
&& \left.-4\sum_{k=0}^{\infty}\left(\begin{array}{c}
3-\alpha \\
k
\end{array}
\right)\left(\frac{-1}{m}\right)^k+
\sum_{k=0}^{\infty}\left(\begin{array}{c}
3-\alpha \\
k
\end{array}
\right)\left(\frac{-2}{m}\right)^k
\right]
{\nonumber}
\end{eqnarray}
leading to
\begin{eqnarray}
q_{j,j-m}   & = &m^{3-\alpha}
\left[\sum_{k=4}^{\infty}\left(\begin{array}{c}
3-\alpha \\
k
\end{array}
\right)\left(\frac{2}{m}\right)^k
-4\sum_{k=4}^{\infty}\left(\begin{array}{c}
3-\alpha \\
k
\end{array}
\right)\left(\frac{1}{m}\right)^k\right.{\nonumber}\\
&& \left.-4\sum_{k=4}^{\infty}\left(\begin{array}{c}
3-\alpha \\
k
\end{array}
\right)\left(\frac{-1}{m}\right)^k+
\sum_{k=4}^{\infty}\left(\begin{array}{c}
3-\alpha \\
k
\end{array}
\right)\left(\frac{-2}{m}\right)^k
\right] {\nonumber}\\
& = & m^{3-\alpha}
\left[ \frac{(3-\alpha)(3-\alpha-1)(3-\alpha-2)(3-\alpha-3)}{4!}\frac{24}{m^4}
  + \dots
\right]
{\nonumber} \\
& = & \frac{1}{m^{\alpha-1}}
\left[ \frac{(3-\alpha)(2-\alpha)(1-\alpha)(-\alpha)}{4!}\frac{24}{m^2}
  + \dots
\right].
\label{ba}
\end{eqnarray}
Considering (\ref{ba}) and noting that the $k$ odd terms of the series cancel,
the properties (a)
can be easily obtained.

(b) In order to compute the series, let us first compute the sum of the first
$M-1$ terms. We have
$$\ds{\sum_{m=2}^{M}} q_{j,j-m}= -3+3\times 2^{3-\alpha}-3^{3-\alpha} + s_M,$$
where
$$
s_M = -(M-1)^{3-\alpha}+3M^{3-\alpha}-3(M+1)^{3-\alpha}+(M+2)^{3-\alpha}.
$$
Similar to what is done in (a) we can write
\begin{eqnarray}
s_M  & = & M^{3-\alpha}
\left[\left(1+\frac{2}{M}\right)^{3-\alpha}-3\left(1+\frac{1}{M}\right)^{3-\alpha}
+3-\left(1-\frac{1}{M}\right)^{3-\alpha}
\right]{\nonumber}\\
&= &M^{3-\alpha}
\left[\sum_{k=0}^{\infty}\left(\begin{array}{c}
3-\alpha \\
k
\end{array}
\right)\left(\frac{2}{M}\right)^k
-3\sum_{k=0}^{\infty}\left(\begin{array}{c}
3-\alpha \\
k
\end{array}
\right)\left(\frac{1}{M}\right)^k+3\right.{\nonumber}\\
&& \left.-\sum_{k=0}^{\infty}\left(\begin{array}{c}
3-\alpha \\
k
\end{array}
\right)\left(\frac{-1}{M}\right)^k
\right].
{\nonumber}
\end{eqnarray}
Therefore
\begin{eqnarray}
s_M   & = &M^{3-\alpha}
\left[\sum_{k=3}^{\infty}\left(\begin{array}{c}
3-\alpha \\
k
\end{array}
\right)\left(\frac{2}{M}\right)^k
-3\sum_{k=3}^{\infty}\left(\begin{array}{c}
3-\alpha \\
k
\end{array}
\right)\left(\frac{1}{M}\right)^k\right.{\nonumber}\\
&& \left.-\sum_{k=3}^{\infty}\left(\begin{array}{c}
3-\alpha \\
k
\end{array}
\right)\left(\frac{-1}{M}\right)^k
\right] {\nonumber}\\
& = & M^{3-\alpha}
\left[ \frac{(3-\alpha)(2-\alpha)(1-\alpha)}{3!}\frac{6}{M^3}
  + \dots
\right]{\nonumber}\\
& = &\frac{1}{M^{\alpha-1}}
\left[ \frac{(3-\alpha)(2-\alpha)(1-\alpha)}{3!}\frac{6}{M}
  + \dots
\right].
\label{sm}
\end{eqnarray}
Clearly, we can conclude that $\lim_{M\rightarrow \infty}s_M=0$.
Hence,
$$
\sum_{m=2}^{\infty}q_{j,j-m} = \lim_{M\rightarrow \infty}
\sum_{m=2}^{M}q_{j,j-m}
= -3+3\times 2^{3-\alpha}-3^{3-\alpha}. 
$$

(c) This result comes immediately  from (b) and from the fact that
$q_{j,j+1}+q_{j,j} + q_{j,j-1} = 3 - 3 \times 2^{3-\alpha}+3^{3-\alpha}$.

\

{\bf Remark}: Note that, the previous result lead us to conclude
the series, defining the operator (\ref{intsumqjm}), converges absolutely when we
have a bounded function $u$.

\

The next theorem states the method is unconditionally stable for $\tau \in [1/2,1]$.

\newtheorem{theorem5}[theorem1]{Theorem}
\begin{theorem5}
The weighted numerical method (\ref{cn}) is unconditionally von Neumann stable for $\tau \in [1/2,1]$.
\end{theorem5}

{\bf Proof}:
 Let us insert the mode $\kappa^n {\rm e}^{ij\theta}$ into (\ref{errorcn}).
We obtain the following
\begin{eqnarray*}
&& \kappa^{n+1}(\theta)
\left[{\rm e}^{ij\theta}-\tau\frac{\mu^\alpha}{\Gamma(4-\alpha)}\sum_{m=-1}^{\infty} q_{j,j-m}
  {\rm e}^{i(j-m)\theta}\right]\\
&&=
\kappa^{n}(\theta)\left[{\rm e}^{ij\theta}+(1-\tau)\frac{\mu^\alpha}
{\Gamma(4-\alpha)}\sum_{m=-1}^{\infty} q_{j,j-m}
  {\rm e}^{i(j-m)\theta}\right].
\end{eqnarray*}
The amplification factor is given by
$$
\kappa(\theta)
\left[1-\tau\frac{\mu^\alpha}{\Gamma(4-\alpha)}\sum_{m=-1}^{\infty} q_{j,j-m}
  {\rm e}^{-im\theta}\right]
=
\left[1+(1-\tau)\frac{\mu^\alpha}{\Gamma(4-\alpha)}\sum_{m=-1}^{\infty} q_{j,j-m}
  {\rm e}^{-im\theta}\right].
$$
Therefore $|\kappa(\theta)| \leq 1$ if and only if the real part of the
series is negative, that is,
$$
\sum_{m=-1}^{\infty} q_{j,j-m} \cos(m\theta) \leq 0,
$$
since the imaginary part of the right side is smaller
for  $\tau \in [1/2,1]$, because $\tau \geq 1-\tau$.
We can write
\begin{eqnarray}
\sum_{m=-1}^{\infty} q_{j,j-m} \cos(m\theta)& = & (q_{j,j+1}+q_{j,j-1})
\cos(\theta) + q_{j,j} {\nonumber}\\
&& + \sum_{m=2}^{\infty} q_{j,j-m} \cos(m\theta).
\end{eqnarray}
Since $q_{j,j+1}+q_{j,j-1} \geq 0, $ and  $ q_{j,j-m} \geq 0$ for $m\geq 2$,
\begin{eqnarray}
\sum_{m=-1}^{\infty} q_{j,j-m} \cos(m\theta) & \leq & (q_{j,j+1}+q_{j,j-1}) +
q_{j,j}+\sum_{m=2}^{\infty} q_{j,j-m} .
\end{eqnarray}
Now using Lemma 3. (c), we obtain
\begin{eqnarray}
\sum_{m=-1}^{\infty} q_{j,j-m} \cos(m\theta)
& \leq  & 0 .
\end{eqnarray}
$\square$

\section{Matricial form}

We start to describe the matricial form of the numerical method,
taking in consideration that to implement the numerical method we need
to have a computational bounded domain.  Let us assume we consider the
computational domain $[a,b]$, where the mesh is defined as
$x_j=a+j\Delta x$ and we assume we have
$$
u(a,t) = 0, \qquad \mbox{and} \qquad u(b,t) = g_b(t) \qquad \mbox{given} 
$$

It is straightforward to conclude, that if $u(a,t) =0$, the problem is
equivalent to a problem defined in the whole real line with the solution 
zero for $x \leq a$.

The numerical method can be written in the matricial form
\begin{eqnarray}
\left(I-\tau\frac{\mu^\alpha}{\Gamma(4-\alpha)}Q\right){\bf U}^{n+1} & = & 
\left(I+(1-\tau)\frac{\mu^\alpha}{\Gamma(4-\alpha)}Q\right){\bf U}^{n} {\nonumber}\\
&&+ \frac{\mu^\alpha}{\Gamma(4-\alpha)}\left(\tau{\bf b}^{n+1}+(1-\tau){\bf b}^{n}\right)
+{\bf p}^{n+\tau},
\label{mf}
\end{eqnarray}
where
$\qquad {\bf p}^{n+\tau}= \left[ \Delta t \tau p_1^{n+1}+(1-\tau)p_1^{n} \dots \Delta t  \tau p_{N-1}^{n+1}+(1-\tau)p_{N-1}^n\right]^T$, 
 \newline ${\bf U}^{n}=\left[U_1^n \dots U_{N-1}^n\right]^T$,
${\bf b}^{n}$ contains the boundary values, $\mu^\alpha$ is a diagonal matrix
with entries $\mu_j^\alpha$ and $Q$ is related to the
fractional operator. The matrix $Q=[Q_{j,k}]$ has the following structure
$$
Q_{j,k} =\left\{
\begin{array}{ll}
q_{j,k}, & 1 \leq k \leq j-1 \\
q_{j,j},  & k=j\\
q_{j,j+1}, & k=j+1\\
0, & k>j+1. 
\end{array}
\right.
$$
Finally the vector ${\bf b}^{n}$ is given by
$$
{b}_j^n =\left\{
\begin{array}{ll}
0, & j=1, \dots, N-2 \\
q_{j,j+1}U_N^n, & j=N-1.
\end{array}
\right.
$$
assuming that $U_0^n=0$ and $U_N^n=g_b(t_n)$.

{\bf Remark}:
From Lemma 4,  for $q_{j,j-1}\geq 0$ (i.e. 
$\alpha > 1.5545 $), we can also easily prove our numerical method is
unconditionally stable by the Gerschgorin's theorem applied to the iterative matrix.

\section{Numerical implementation}

The numerical experiments are carried out in two parts. First, we
verify the accuracy and order of convergence of the numerical method to confirm 
the theoreticall results presented in the previous sections. Then
a physical application is considered to reveal some of the physical
phenomena, from anomalous to mormal diffusion.

Consider the vectors $U_{app} (\Delta x) = (U_0, \dots, U_N)$, where
$U_j$ is the approximate solution, for $x_j= x_0+j\Delta x$, $j=0,
\dots, N$ at a certain time $t$, and $u_{ex}(\Delta x)= (u(x_0,t),
\dots, u(x_N,t))$, where $u$ is the exact solution.  The error is
defined by the $l_\infty$ norm as,
\begin{equation}
||u_{ex}(\Delta x)-U_{app}(\Delta x)||_\infty = \max_{0 \leq j \leq N} 
\left| u(x_j,t)-U_j\right|.
\label{ie}
\end{equation}

{\bf Example 1.}
Consider the problem with initial condition $u(x,0)= 4x^2(2-x)^2$, $0<x<2$ and
zero otherwise.
Let 
\begin{equation}
d(x)=\frac{1}{4}\Gamma(5-\alpha)x^{\alpha},
\label{dp1}
\end{equation}
and
\begin{equation}
p(x,t) = -4{\rm e}^{-t}x^2\left[7(2-x)^2+2\alpha(\alpha-7)+6x\alpha\right].
\label{pp1}
\end{equation}
The exact solution is given by
$
u(x,t) = 4 {\rm e}^{-t}x^2(2-x)^2$, for $ 0 \leq x \leq 2$,
and zero otherwise.

In Table 1, we show the behaviour of the error (\ref{ie}) for different values of $\tau$
and for $\Delta t=\Delta x = 1/30$ for the problem (\ref{dp1})--(\ref{pp1}).

\begin{table}[h]
  \begin{center}
    \begin{tabular}{ccccc}\hline
         $\tau$ & $\alpha=1.2$&  $\alpha=1.4$ &  $\alpha=1.5$ &  $\alpha=1.8$ \\\hline
         0.5&\textbf{\small4.0277$\times10^{-3}$}  &\textbf{\small3.4191$\times10^{-3}$}
         &\textbf{\small3.1944$\times10^{-3}$}&\textbf{\small2.4542$\times10^{-3}$}  \\
         0.6&5.6194$\times10^{-3}$  &4.9682$\times10^{-3}$  &4.7877$\times10^{-3}$ &4.2856$\times10^{-3}$ \\
         0.7&7.5094$\times10^{-3}$  &6.7573$\times10^{-3}$  &6.5920$\times10^{-3}$ &6.2510$\times10^{-3}$ \\
         0.8&9.5429$\times10^{-3}$  &8.6634$\times10^{-3}$  &8.4903$\times10^{-3}$ &8.2598$\times10^{-3}$ \\
         0.9&1.1656$\times10^{-2}$  &1.0625$\times10^{-2}$  &1.0435$\times10^{-2}$ &1.0283$\times10^{-2}$\\
         1.0&1.3814$\times10^{-2}$  &1.2615$\times10^{-2}$  &1.2403$\times10^{-2}$ &1.2318$\times10^{-2}$\\\hline
\\
    \end{tabular}\label{table1}
  \end{center}
\caption{Global $ l_{\infty} $ error (\ref{ie}) of time converged solution at \ $t=1$
\  for \ $ \alpha=1.2$, \ $\alpha=1.4$, \ $\alpha=1.5$, \ $\alpha=1.8$ and $\Delta t=\Delta
    x=1/30$.}\vspace{3pt}
\end{table}

The most accurate result is for $\tau=1/2$.
For the same problem,
we observe in Table 2 and Table 3 that for all values of $\alpha$
we have second order convergence as expected, when $\tau=1/2$.

\begin{table}[h]
  \begin{center}
    \begin{tabular}{ccccc}\hline
         $\Delta x$ & $\alpha=1.2$&  Rate &  $\alpha=1.4$ &  Rate \\\hline
         1/5 &1.5310$\times10^{-1}$  &        &1.1950$\times10^{-1}$&  \\
         1/10&3.6239$\times10^{-2}$  &2.0789  &3.0270$\times10^{-2}$&1.9811  \\
         1/20&9.0506$\times10^{-3}$  &2.0015  &7.6627$\times10^{-3}$&1.9820 \\
         1/40&2.2669$\times10^{-3}$  &1.9973  &1.9289$\times10^{-3}$&1.9901 \\ \hline
    \end{tabular}\label{table2}
  \end{center}
    \caption{Global $l_{\infty}$ error (\ref{ie})of time converged solution for four mesh resolutions at $t=1$
    for $\alpha=1.2,\alpha=1.4$, $\Delta t=\Delta
    x$ and $\tau=1/2$.}\vspace{2pt}
\end{table}

\begin{table}[h]
  \begin{center}
    \begin{tabular}{ccccc}\hline
         $\tau$ & $\alpha=1.5$&  Rate &  $\alpha=1.8$ &  Rate \\\hline
         1/5 &1.0884$\times10^{-1}$  &        &7.9651$\times10^{-2}$ & \\
         1/10&2.8101$\times10^{-2}$  &1.9535  &2.0820$\times10^{-2}$ &1.9357 \\
         1/20&7.1358$\times10^{-3}$  &1.9775  &5.4174$\times10^{-3}$ &1.9423 \\
         1/40&1.8050$\times10^{-3}$  &1.9831  &1.3974$\times10^{-3}$ &1.9549 \\\hline
    \end{tabular}\label{table3}
  \end{center}
   \caption{Global $l_{\infty}$ error (\ref{ie}) of time converged solution for four mesh resolutions at $t=1$
    for $\alpha=1.5,\alpha=1.8$, $\Delta t=\Delta
    x$ and $\tau=1/2$.}\vspace{2pt}
\end{table}

{\bf Example 2.} Consider now a second problem with initial condition $u(x,0)= x^\lambda$, $0\leq x\leq1$ and
boundary conditions $u(0,t)=0$ and $u(1,t)={\rm e}^{-t}$.
Let
\begin{equation}
d(x)=\frac{\Gamma(\lambda+1-\alpha)}{\Gamma(\lambda+1)}x^{\alpha+1} \quad \mbox{and} \quad p(x,t) = -(1+x){\rm e}^{-t}x^\lambda.
\label{p2}
\end{equation}
The exact solution of the problem is of the form
\begin{equation}
u(x,t) = {\rm e}^{-t}x^\lambda, \qquad x \in [0,1].
\label{es}
\end{equation}
Although this problem is not defined in the whole real line we have $u(0,t)=0$,
and this can be seen as a problem for which the solution is zero when $x\leq 0$.

In Table 4, we show
the behavior of the error (\ref{ie}) for
different weighted coefficients
$\tau$. We
observe the most accurate behaviour is again for $\tau=1/2$.


\begin{table}[h]
  \begin{center}
    \begin{tabular}{ccccc} \hline
         $\tau$ & $\alpha=1.2$&  $\alpha=1.4$ &  $\alpha=1.5$ &  $\alpha=1.8$ \\\hline
         0.5&\textbf{\small6.4792$\times10^{-5}$}  &\textbf{\small2.9402$\times10^{-5}$}
         &\textbf{\small1.7850$\times10^{-5}$}&\textbf{\small4.0509$\times10^{-6}$}  \\
         0.6&9.6854$\times10^{-4}$  &7.0639$\times10^{-4}$  &6.2104$\times10^{-4}$ &4.5122$\times10^{-4}$ \\
         0.7&1.8609$\times10^{-3}$  &1.3815$\times10^{-3}$  &1.2233$\times10^{-3}$ &9.0545$\times10^{-4}$ \\
         0.8&2.7426$\times10^{-3}$  &2.0533$\times10^{-3}$  &1.8233$\times10^{-3}$ &1.3587$\times10^{-3}$\\
         0.9&3.6143$\times10^{-3}$  &2.7219$\times10^{-3}$  &2.4211$\times10^{-3}$ &1.8110$\times10^{-3}$ \\
         1.0&4.4769$\times10^{-3}$  &3.3870$\times10^{-3}$  &3.0166$\times10^{-3}$ &2.2624$\times10^{-3}$
 \\\hline
\\
    \end{tabular}
  \end{center}
\caption{Global $l_\infty$ error (\ref{ie}) of time converged solution
for the problem (\ref{p2}) calculated by weighted numerical scheme with $\Delta t=\Delta x=1/30,\lambda=3, 0 \leq x \leq 1$
     for different values of $\alpha$ and $\tau$.}\vspace{2pt}
\label{table4}
\end{table}

In Table \ref{table5} we present a comparison between our method and
the methods presented in \cite{tad2006} with the same space and time
steps. The second column shows the absolute value of the largest error
calculated by the Crank-Nicolson scheme (before extrapolation)
presented in \cite{tad2006} at time $t=1.0$ which consists of assuming
the fractional derivative is approximated by the shifted
Gr\"{u}nwald-Letnikov formula. The third column shows the error calculated
by the Crank-Nicolson scheme after a Richardson's extrapolation
presented in \cite{tad2006}. The fourth column shows the largest
absolute error for our numerical scheme with $\tau=0.5$.  Note that
our numerical results are  more accurate than the method
given in \cite{tad2006}.


\begin{table}[h]
  \begin{center}
    \begin{tabular}{cccc} \hline
                $\Delta x$ & CN-GL \cite{tad2006}& Extrapolated CN-GL \cite{tad2006}&Weighted ($\tau=0.5$)\\
                               \hline
                                1/10 & 1.82265$\times10^{-3}$  & 1.77324$\times10^{-4}$  & 3.5504$\times10^{-5}$  \\
                               1/15 & 1.16803$\times10^{-3}$  & 7.85366$\times10^{-5}$  & 1.6197$\times10^{-5}$  \\
                                1/20 & 8.64485$\times10^{-4}$  & 4.40627$\times10^{-5}$  & 9.1072$\times10^{-6}$ \\
                              1/25 & 6.84895$\times10^{-4}$  & 2.82750$\times10^{-5}$  & 5.8030$\times10^{-6}$ \\\hline
\\
    \end{tabular}
  \end{center}
    \caption{Global $l_{\infty}$ error (\ref{ie}) of time converged solution for the second problem
    calculated at $t=1$ for the second problem with $\Delta t=\Delta
    x,\lambda=3,0\leq x\leq1$  and $\alpha=1.8$.}
\label{table5}
\end{table}

To conclude this example we observe the rate of convergence 
of the numerical method for different values of $\tau \neq 1/2$.
The expected convergence rate for $\tau \neq 1/2$ according to section 3 is
$O(\Delta t+\Delta x^2)$.
We consider $\Delta t=\Delta x^2$ to get second order convergence as we observe in
Table 6.

\begin{table}[h]
  \begin{center}
    \begin{tabular}{cccccc} \hline
                                   &$\Delta t$ &$\Delta x$& $\alpha=1.8$ &Rate\\\hline
                                   & 1/25  &1/5  & 7.9325$\times10^{-4}$ &  -   \\
\raisebox{2.3ex}[0pt]{$\tau=0.6$}  & 1/100 &1/10 & 2.1501$\times10^{-4}$ & 1.8834 \\
                                   & 1/400 &1/20 & 5.3710$\times10^{-5}$ & 2.0011 \\
                                   & 1/1600&1/40 & 1.3512$\times10^{-5}$ & 1.9909 \\
           \hline
                                   & 1/25  &1/5 & 1.2837$\times10^{-3}$ & -     \\
\raisebox{2.3ex}[0pt]{$\tau=0.7$}  & 1/100 &1/10& 3.5212$\times10^{-4}$ & 1.8662 \\
                                   & 1/400 &1/20& 8.7892$\times10^{-5}$ & 2.0023 \\
                                   & 1/1600&1/40& 2.2053$\times10^{-5}$ & 1.9948 \\
           \hline
                                   & 1/25  &1/5 & 1.7723$\times10^{-3}$ & -   \\
\raisebox{2.3ex}[0pt]{$\tau=0.8$}  & 1/100 &1/10& 4.8915$\times10^{-4}$ & 1.8573 \\
                                   & 1/400 &1/20& 1.2207$\times10^{-4}$ & 2.0026  \\
                                   & 1/1600&1/40& 3.0594$\times10^{-5}$ & 1.9964 \\
           \hline
                                   & 1/25  &1/5 & 2.2590$\times10^{-3}$ & -   \\
\raisebox{2.3ex}[0pt]{$\tau=0.9$}  & 1/100 &1/10& 6.2608$\times10^{-4}$ & 1.8513 \\
                                   & 1/400 &1/20& 1.5624$\times10^{-4}$ & 2.0026  \\
                                   & 1/1600&1/40& 3.9134$\times10^{-5}$ & 1.9973 \\
           \hline
                                   & 1/25  &1/5  & 2.7438$\times10^{-3}$ & -   \\
\raisebox{2.3ex}[0pt]{$\tau=1.0$}  & 1/100 &1/10 & 7.6292$\times10^{-4}$ & 1.8466 \\
                                   & 1/400 &1/20 & 1.9041$\times10^{-4}$ & 2.0024 \\
                                   & 1/1600&1/40 & 4.7674$\times10^{-5}$ & 1.9978 \\ \hline
    \end{tabular}
  \end{center}
    \caption{Global $l_{\infty}$ error (\ref{ie}) of time converged solution for the
    numerical scheme (\ref{cn}) at $t=1$, for $\Delta t= \Delta x^2$
and $\alpha=1.8$, $\lambda=3$ and different values of $\tau$.}\vspace{3pt}
\label{table6}
\end{table}

\begin{figure}[h]
\centerline{
\psfig{figure=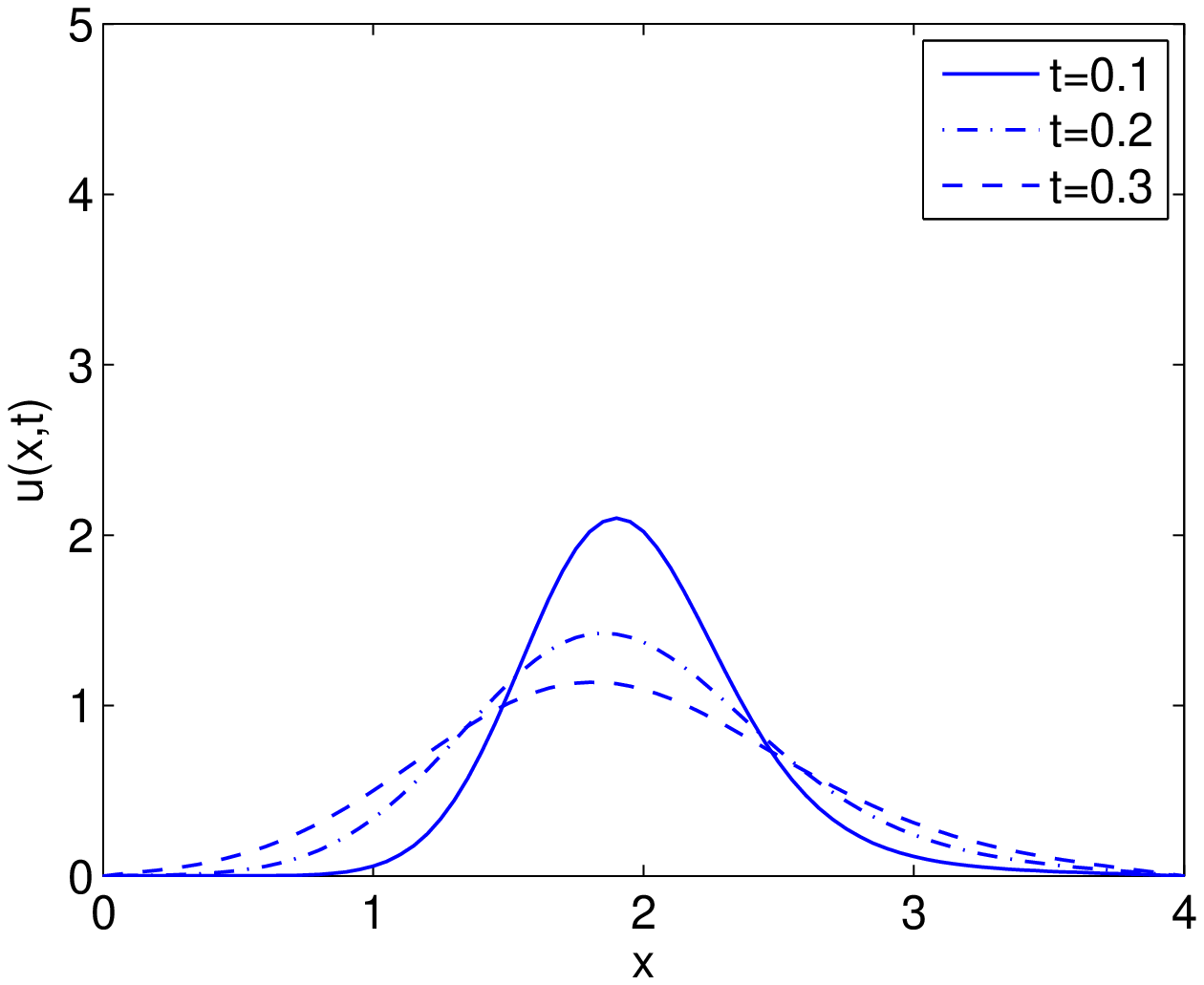,height=6.5cm,width=7.5cm}
\psfig{figure=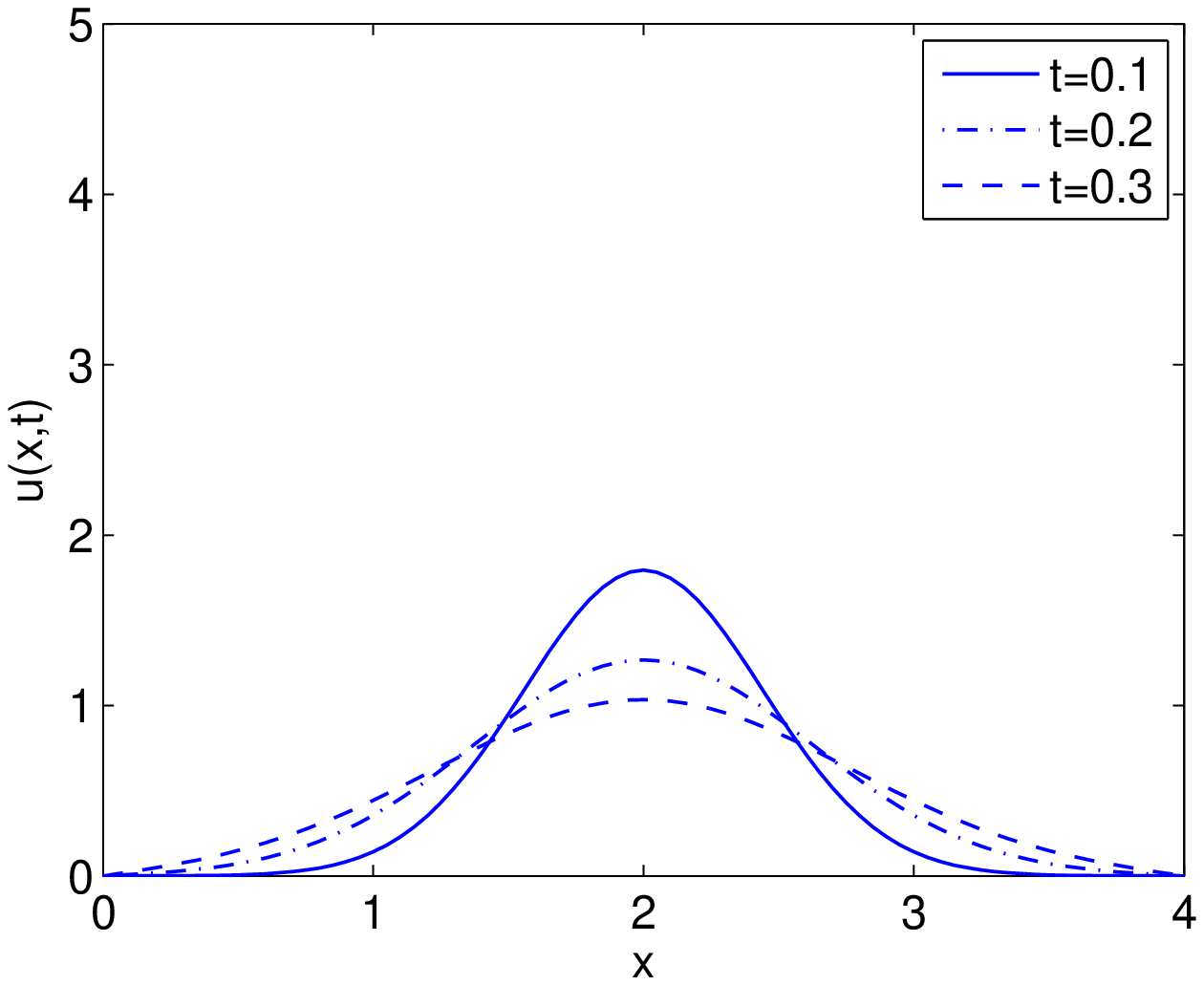,height=6.5cm,width=7.5cm}
}
\hspace*{2.5cm} (c) $\alpha=1.8$\hspace*{5cm} 
(d) $\alpha=1.999$
\caption{The evolution of $u(x,t)$
for different anomalous diffusion coefficients $\alpha$ at different times.}
\label{figfinal}
\end{figure}

{\bf Example 3.} Finally,
in order to reveal the dynamics behavior of the diffusion equation
(\ref{fde}), in this example we consider equation (\ref{fde}) without the source
function (which means $p(x,t)=0$) on a finite domain $[0,4]$.  We consider
the Gaussian function
\begin{displaymath}
u(x,0)=\frac{1}{\sigma\sqrt{2\pi}}\exp(-\frac{(x-2)^2}{2\sigma^2})
\end{displaymath}
as the initial condition,
the diffusion coefficient $d(x)=1$  and the
boundary conditions
$
 u(0,t) = u(4,t) = 0.
$
The numerical results for this example are calculated
by the weighted scheme with $\tau=1/2$.  In this test, we take
$\sigma=0.01$. The evolution of the non-Fickian diffusion processes
for different values of $\alpha$  are given in  Fig \ref{figfinal}. The
anomalous diffusion parameter exhibits the extent of the long tail
diffusion processes of problem (\ref{fde}). The non-Fickian behavior
gradually disappear when $\alpha\rightarrow 2$. This is
consistent with the experimental results \cite{ben2000a, hua2006, pac2000,zho2003}. Again the
validity of our numerical methods is confirmed.

\begin{figure}[h]
\centerline{
\psfig{figure=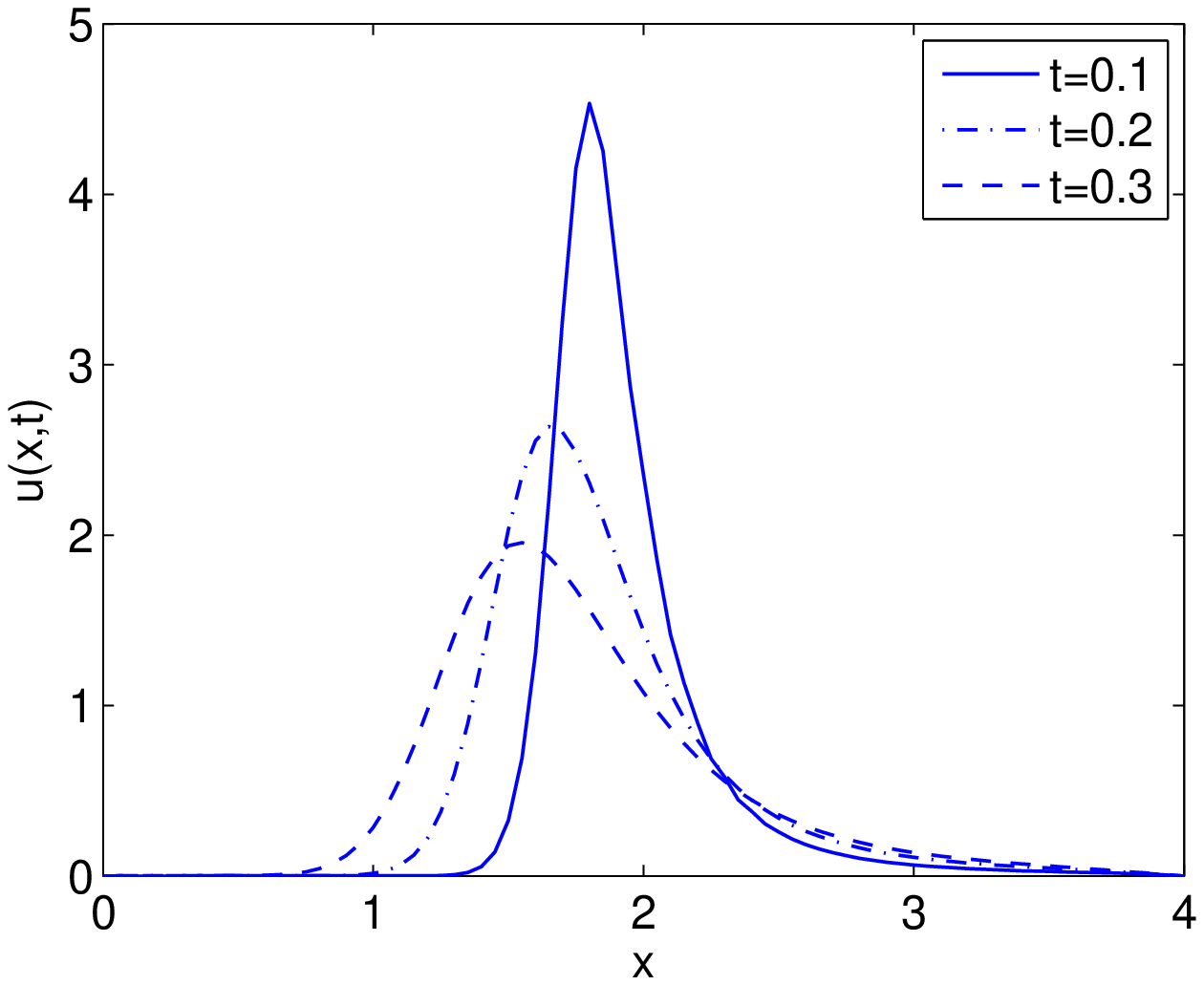,height=6.5cm,width=7.5cm}
\psfig{figure=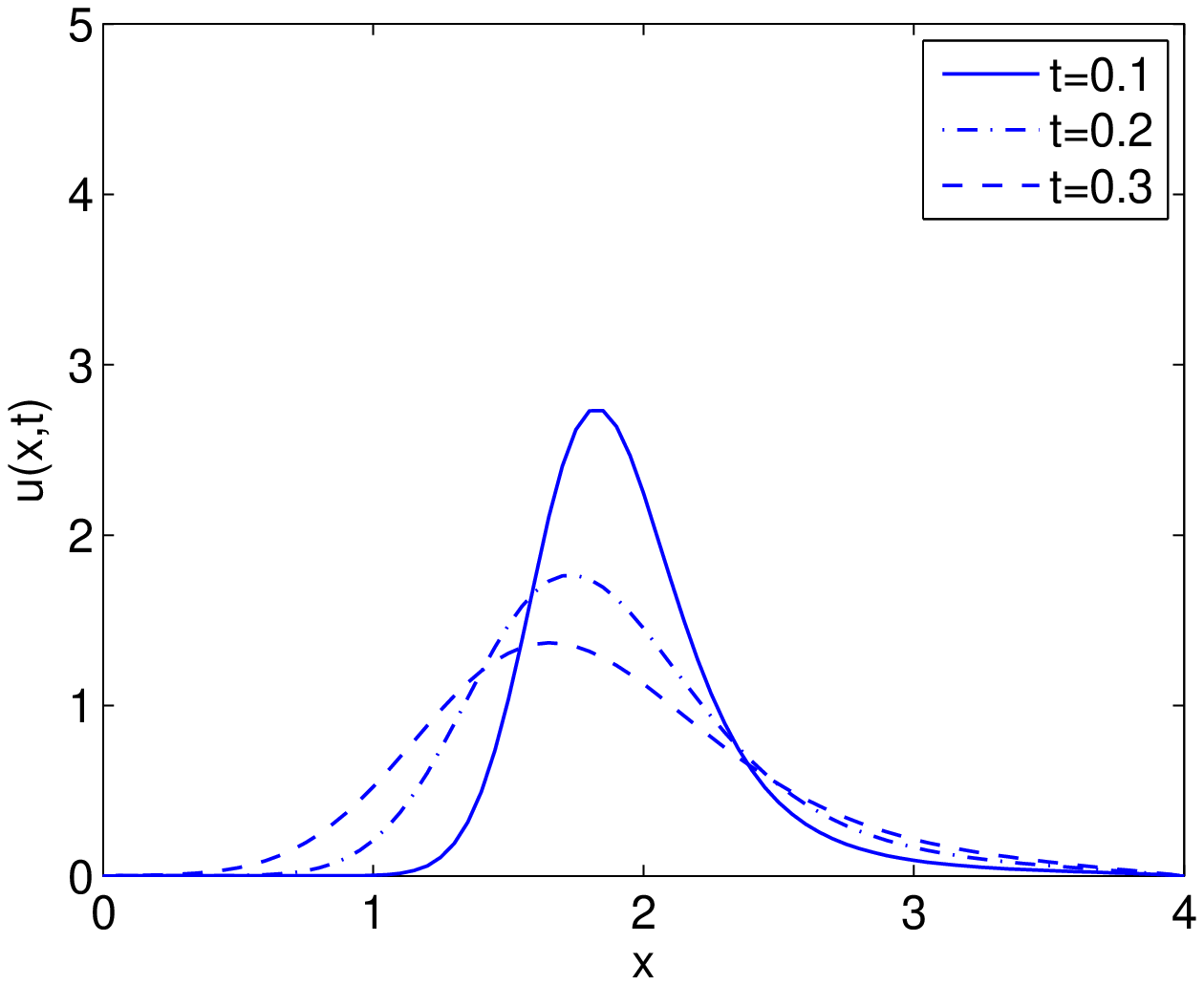,height=6.5cm,width=7.5cm}
}
\hspace*{2.5cm}(a) $\alpha=1.4$\hspace*{5cm} 
(b) $\alpha=1.6$
\end{figure}

\section{Conclusions}

We have derived a weighted numerical method for the fractional diffusion equation
based on the Riemann-Liouville derivative defined  in an unbounded domain.
The numerical method is second order accurate for $\tau=1/2$
and first order accurate for $\tau \in (1/2,1]$ because of the time discretization.
We have proved theoretically the method converges by showing consistency
and von Neumann stability. In the end we have presented test problems which are in
agreement with the theoretical results.

\end{document}